# A reformulation to Embedding a Neural Network in a linear program without integer variables.

Héctor G.-de-Alba [1], Andres Tellez [2], Cirpiano Santos [3], Emmanuel Gómez [4]

[1,2,3,4] Tecnologico de Monterrey, Escuela de Ingeniería y Ciencias, Ave. Eugenio Garza Sada 2501, Monterrey, 64849, Nuevo Léon, México

[1]hrgarcia@tec.mx, [2]andrestellez84@tec.mx, [3]cipriano.santos@tec.mx, [4]emmanuel.gr@tec.mx

**Abstract**

In this technical report, a new formulation for embedding a neural network into an optimization model is described. This formulation does not require binary variables to properly compute the output of the neural network for specific types of problems. Preliminary experiments show that this reformulation resulted in faster computation times when solving a proposed showcase model, in which non-linearity is necessary to be computed. This is in comparison with the classic formulation and off-the-shelf tools of commercial solvers.

It is worth noting that this work is still in progress, and further analysis is needed to properly define the benefits and limitations of the proposed formulation.

**Introduction**

Mathematical programming is a widely used optimization technique that has been employed for several decades to address complex combinatorial problems [1]. Mathematical programming models describe systems using linear expressions that define a space of solutions feasible to the problem at hand (referred to as constraints) and then explore said space to identify the optimal value of a given function (referred to as the objective function). The decision variables of these problems can be non-negative, continuous, and/or integer. The inclusion of integer variables in an optimization model results in a mixed-integer linear program (MILP), which is generally considered to have NP-hard complexity [2],[3].

Additionally, in order to find a global optimal solution when solving a mathematical optimization model, certain criteria must be met beforehand, one of which is that the solution space is convex. Linear expressions are always preferred to ensure convexity. The motivation for this work is to facilitate the linearization of non-linear expressions by embedding neural networks (NN) that have been trained to emulate said functions [4] into optimization models, and to do this without the need for additional integer variables. Work related to the use of neural networks in optimization models can be viewed in [5], [6], and [7].

In this technical report, we first discuss basic concepts regarding neural networks, then recapitulate the existing methods in the literature to embed a neural network. Afterward, we present our proposed method to embed neural networks and explain its implications. Lastly, we demonstrate an example in which the proposed neural network embedding is applied and compared to the classical method and other off-the-shelf commercial solutions in preliminary experiments.

It is worth clarifying that an appropriate degree of knowledge in the fields of optimization, linear programming, and machine learning is expected from the reader, as this technical report is aimed at such readers.

**Neural network Basic Concepts**

An artificial neural network (NN) is a machine learning model inspired by the functioning of human brain cells and how they interconnect to process information. A NN consists of various nodes (analogous to neurons), which are units that perform computations by receiving a set of inputs and processing them through a specific function to produce an output. Neurons are organized into layers: an input layer that receives values from the user, one or more hidden layers to process the information (working as a black box to the user), and an output layer that presents the resulting values of the computations done by the NN to the user. In this arrangement, the output produced by the neurons of a specific layer becomes the input for the nodes in the following layer, and this communication between nodes is facilitated through connections (analogous to axons).

Considering this design, input data is fed into the input layer and undergoes a series of transformations as it passes through the hidden layers until the network produces an output. In this regard, neural networks are highly versatile and can be used to perform a wide range of tasks. Figure 1 in the diagram represents a visualization of a NN.

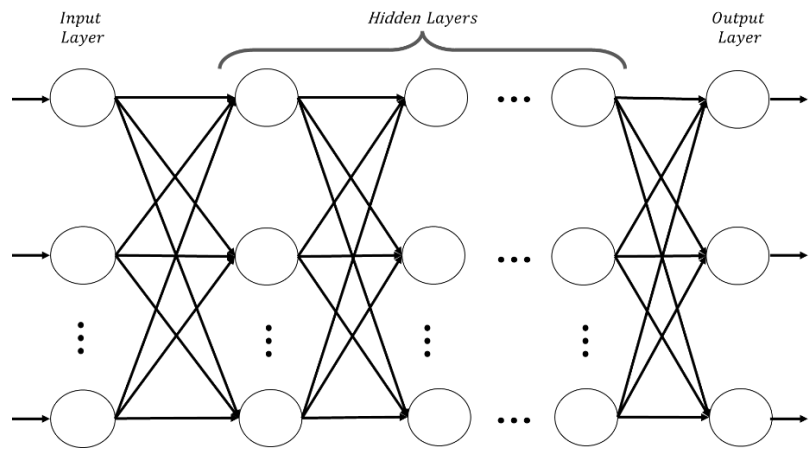

*Figure 1: Diagram of a NN*

In the diagram in Figure 1, it can be observed that all nodes from one layer exclusively send information to all nodes in the immediately following layer. This architecture/arrangement of nodes and connections is referred to as a "feed-forward densely connected NN." It should be clarified that other architectures exist, but this particular one is explored for the application intended for the neural network embeddings proposed in this document.

Additionally, determining the number of neurons and nodes to be used in the network is a decision process typically undertaken in an artisanal fashion, as there is no way to beforehand know the optimal NN architecture for a specific task.

As mentioned earlier, the information outputted by a node will be broadcasted to the nodes in the next layer through connections. This information may be intentionally "modulated" by these connections through the use of weights. In other words, when the output of a node is transmitted through the connections to the nodes in the next layers, the weights of the connections will multiply the broadcasted value. This can be observed in the diagram in Figure 2.

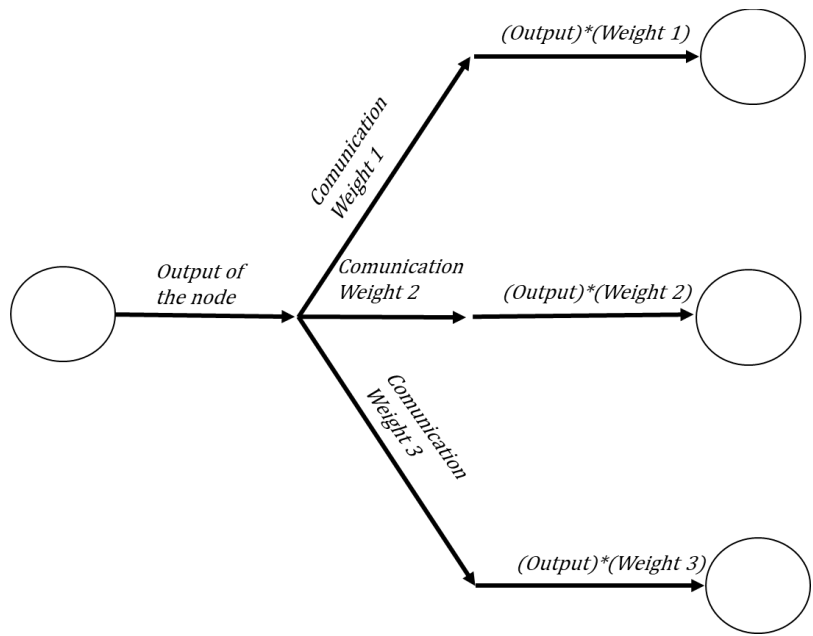

Figure 2: weights diagram

Concerning how nodes compute their output value, it is achieved by processing the sum of the weighted inputs that the node receives. The function used for this process is known as an activation function. There are various types of functions that can be considered for use as an activation function; examples include the sigmoid function, hyperbolic function, among many others. For this document, we will consider the rectifier linear unit (ReLU) function, defined as $f(x) = \max(0, x)$. Additionally, to regulate the responsiveness of neurons to their inputs, each of them will have a bias value added to the weighted values before being processed by the activation function. These aspects are illustrated in the diagram in Figure 3.

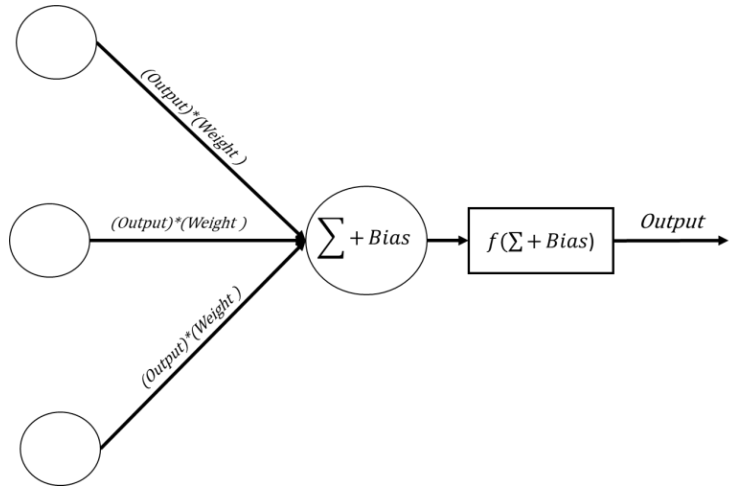

*Figure 3: node activtion diagram*

The weight values of the connections and the bias values of each neuron are adjusted during the learning process to enhance the network's performance in the designated task. This process, known as training, involves presenting a set of inputs along with their corresponding expected outputs to the network. Subsequently, an optimization algorithm, typically the backpropagation algorithm, is applied to appropriately adjust the weights and biases, aiming to minimize the errors between the produced outputs and the expected outputs.

Once the training reaches an adequate degree, characterized by achieving a low error, the neural network captures complex patterns in its internal layers to establish the proper relationship between input values and expected output values. These capabilities allow these models to approximate any function, as asserted by the Universal Approximation Theorem [4], and even learn from seemingly chaotic datasets, provided an appropriate number of layers and nodes are considered for the neural network.

These capabilities are leveraged by embedding a fully trained neural network into an optimization model. For further information on neural networks, refer to [8] and [9].

**Neural network embedding**

In this section, we explain the utilization of neural networks (NN) embedded into optimization models as outlined in current literature. To perform neural network embeddings, the weights of the connections, the biases of the neurons, and the architecture of the network need to be supplied to the optimization model as parameters. All of these parameters should have been acquired through prior training of the network using samples of the inputs that the NN may receive in the optimization model and the corresponding outputs that the network should produce. These parameters are then incorporated into constraints that emulate the behavior of the neural network. The following provides a definition of these parameters:

- $I$: number of hidden layers in the network.

- $J_i$ : number of nodes in layer $i$; with $i = 0 \ldots I+1$, in this regard, layer $i = 0$ represents the input layer and $i = I+1$ the output layer.
- $W_{\hat{j},i,j}$: weight of the connection that goes from the node $\hat{j}$ in the layer $i-1$, to the node $j$ in the layer $i$; with $i = 1 \ldots I+1$, $j = 1 \ldots J_i$ and $\hat{j} = 1 \ldots J_{i-1}$.
- $B_{i,j}$: bias of the neuron $j$ of the layer $i$, with $i = 1 \ldots I, j = 1 \ldots J_i$.

For this document, we will consider a network with a single input and a single output, representing a network that will approximate a function $G(x)$, which can be any type of non-linear function; because of these, we have that $J_0 = J_{I+1} = 1$. A representation of the architecture of this network and the placement of its parameters can be view in the figure 4 diagram:

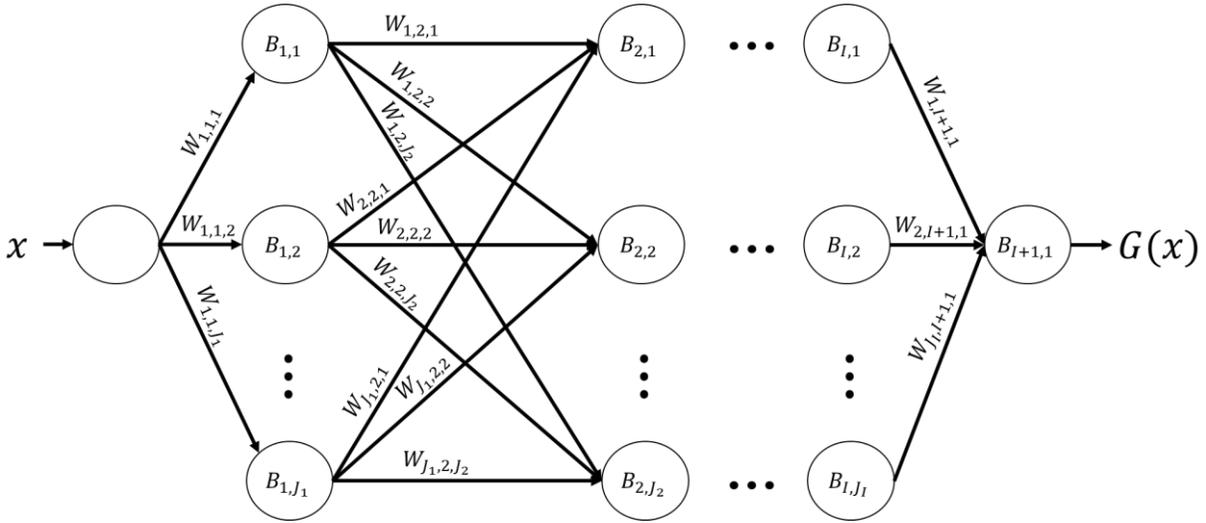

*Figure 4: Representation and placement of parameters for the considered NN*

It is worth noticing that the output node of the network, in which the output $G(x)$ is computed, also has a bias value, even thou in this node the activation function will not be applied, and just the sum of the weighted inputs plus the bias will be computed in order to produce the output of the network; this is a standard practice in the design of neural networks, as if the activation function was applied to the nodes of the output layer, the set of values that the NN can yield as outputs will be limited by said function $f(x)$, which can be undesired behavior if the domain of $f(x)$ doesn't match the domain of $G(x)$.

For the formulation of the constraints that will be used to emulate the behavior of the NN, we need to consider the following decision variables:

- $\sigma_{i,j}$ : output value of the node $j$ of the hidden layer $i$; $\forall\ i = 1 \ldots I, j = 1 \ldots J_i$, and $\sigma_{i,j} \geq 0$.
- $y_{i,j}$: binary variable that takes the value of 1 when $\sigma_{i,j} > 0$, and takes the value of 0 when $\sigma_{i,j} \leq 0$; $\forall\ i = 1 \ldots I, j = 1 \ldots J_i$, and $y_{i,j} \in \{0,1\}$.
- $\theta$: output value of the network; $\theta \in \mathbb{R}$.
- $x$: Input value of the network; $x \in \mathbb{R}$.

It can be noticed that, for each node of the hidden layers of the neural network, its output will be a non-negative value, this is because of the usage of the ReLU function as activation function;

because of this, the embeddings shown in this document can only be used with NN trained with such activation functions. Considering these decision variables, the constraints used to embed the NN into the optimization model will be as follow:

$$\sigma_{1,j} \geq W_{1,1,j}x + B_{1,j} \qquad \forall\, j = 1 \ldots J_i \qquad (1)$$

$$\sigma_{1,j} \leq My_{1,j} \qquad \forall\, j = 1 \ldots J_i \qquad (2)$$

$$\sigma_{1,j} \leq W_{1,1,j}x + B_{1,j} + M(1 - y_{i,j}) \qquad \forall\, j = 1 \ldots J_i \qquad (3)$$

$$\sigma_{i,j} \geq \sum_{\hat{j}=1}^{J_{i-1}} (W_{\hat{j},i,j}\sigma_{i-1,\hat{j}}) + B_{i,j} \qquad \begin{array}{l}\forall i = 2 \ldots I,\\ j = 1 \ldots J_i\end{array} \qquad (4)$$

$$\sigma_{i,j} \leq My_{i,j} \qquad \begin{array}{l}\forall i = 2 \ldots I,\\ j = 1 \ldots J_i\end{array} \qquad (5)$$

$$\sigma_{i,j} \leq \sum_{\hat{j}=1}^{J_{i-1}} (W_{\hat{j},i,j}\sigma_{i-1,\hat{j}}) + B_{i,j} + M(1 - y_{i,j}) \qquad \begin{array}{l}\forall i = 2 \ldots I,\\ j = 1 \ldots J_i\end{array} \qquad (6)$$

$$\theta = \sum_{\hat{j}=1}^{J_I} (W_{\hat{j},I+1,1}\sigma_{I,\hat{j}}) + B_{I+1,1} \qquad (7)$$

Where M is a big constant. Constraints (1) to (3) compute the output of the nodes that are in the first hidden layer, which directly receive the input of the network. Constraints (4) to (6) compute the output of the nodes for the rest of the hidden layers, which receive inputs from other nodes, provided that the network was designed with multiple hidden layers. Lastly, constraint (7) computes the output of the network.

It's worth noticing that the role that the auxiliary binary variables $y_{i,j}$ play in this formulation is to activate and deactivate the constraints in which those variables are present, in order to properly bound each node output value; in this regard, the constraints (2) and (5) bound the output of the node to 0 when the sum of the weighted inputs that receive the node and its bias result in a non-positive value, while the constraints (3) and (6) bound the output of the node to the value resulting of the sum of the weighted inputs that receive the node and its bias when this value is positive. Because of this, the embedding can be applied in any optimization model, but will result in an increase in computation time for said model because of the increase in the number of binary variables, by increasing the ramifications needed in the algorithms used to solve mixed integer models. In order to prevent this increase in computational time, we propose an alternative and more compact formulation that can only be applied in specific scenarios, which will be described in the next section.

**Proposed Neural network embedding.**

The proposed formulation considers the same parameters and variables, except for the auxiliary binary variable, and goes as follows:

$$\sigma_{1,j} \geq W_{1,1,j}x + B_{1,j} \qquad \forall\, j = 1 \ldots J_i \qquad (8)$$

$$\sigma_{i,j} \geq \sum_{\hat{j}=1}^{J_{i-1}} \left(W_{\hat{j},i,j}\sigma_{i-1,\hat{j}}\right) + B_{i,j} \qquad \forall i = 2 \ldots I, \quad j = 1 \ldots J_i \qquad (9)$$

$$\theta = \sum_{\hat{j}=1}^{J_I} \left(W_{\hat{j},I+1,1}\sigma_{I,\hat{j}}\right) + B_{I+1,1} \qquad (10)$$

As observed, this formulation considers some of the original constraints and excludes those related to the auxiliary variables, which are unnecessary for computing the output values of the network. However, in return, its applicability will be limited to models in which the output of the network is intended to be minimized by the model, directly or indirectly. Additionally, the values of the parameters, denoted as $W_{\hat{j},i,j}$, must be non-negative, a condition that must be established during the training of the neural network. Failure to meet these criteria may result in the neural network not being properly bounded, leading to unexpected values.

To illustrate the usage of the proposed embedding, a simple example will be presented in the next section.

**Application example.**

In order to showcase the proposed formulations for embedded a neural network, let's consider a variation of the knapsack problem with $n$ different classes of items, contained by the set $C = \{1 \ldots n\}$. For each class of item $c \in C$, a certain number of copies of it will exist, denoted by the parameter $m_c$, all of them having the same characteristics, such as a size(volume) $s_c$, and a profit value $v_c$. As the typical knapsack, the objective is to maximize the profit value obtained by selecting as many items as possible from each class, which will be limited by the size of the knapsack, denoted by the parameter $S$; for practical reasons, the geometry of the objects will be disregarded when fitting them in the knapsack. In addition, selecting multiple items of the same class will be penalized in a quadratic way, and this penalization might be different for each item class, which will be denoted by the parameter $p_c$. The described problem can be formulated as a non-lineal integer model. The decision variables for the model will be:

- $X_c$ = number of items to be selected from the class $c$; $c = 1 \ldots n$, and $X_c \in \mathbb{N}$.

Considering this, the corresponding model will be:

$$max \sum_{c=1}^{n} (v_c X_c - p_c(X_c^2 - X_c)) \qquad (11)$$

Subject to:

$$\sum_{c=1}^{n} s_c X_c \leq S \qquad (12)$$

$$X_c \leq m_c \qquad \forall c = 1 \ldots n \qquad (13)$$

The expression (11) defines the objective function, it is worth noticing that $p_c(X_c^2 - X_c)$ represents the penalization function incurred by selecting multiple copies of the same item. Expression (12) and (13) are the constraints of the problem, that describe how there is a limit to the total size of the items to select that the knapsack can fit inside, and that there is a limited number of items that can be taken from each class, respectively.

The expression $X_c^2$ is non-linear and might be difficult to tackle for some solvers, and it will be the focus of the NN implementation of the embedding proposed in this document. In order to do this additional decision variables and constraints will need to be incorporated into the model to compute $X_c^2$ by the embedding. Because every single $X_c$ will need to compute its own corresponding $X_c^2$, then an embedding for every single each of them will be necessary, hence, the next variables will be needed for the embedding:

- $\sigma_{i,j,c}$ : output value of the node $j$ of the hidden layer $i$ for the object class $c$; with $i = 1 \ldots I, j = 1 \ldots J_i$, and $\sigma_{i,j} \geq 0$.
- $\theta_c$: output value of the network for the object class $c$; $\theta \in \mathbb{R}$.

In this case, the input values for the embeddings will be the different decision variables $X_c$, and will produce the corresponding $X_c^2$ by the means of the variables $\theta_c$. It is worth noticing that parameters of the connection's weights $W_{j,i,j}$ and biases $B_{i,j}$ will be the same for all the embeddings used to compute every single $X_c^2$. Considering this, the corresponding model will be:

$$max \sum_{c=1}^{n}(v_c X_c - p_c(\theta_c - X_c)) \tag{14}$$

Subject to:

$$\sum_{c=1}^{n} s_c X_c \leq S \tag{15}$$

$$X_c \leq m_c \quad \forall\, c = 1 \ldots n \tag{16}$$

$$\sigma_{1,j,c} \geq W_{1,1,j}X_c + B_{1,j} \quad \begin{array}{l}\forall\, j = 1 \ldots J_i,\\ c = 1 \ldots n\end{array} \tag{17}$$

$$\sigma_{i,j,c} \geq \sum_{\hat{j}=1}^{J_{i-1}}(W_{\hat{j},i,j}\sigma_{i-1,\hat{j},c}) + B_{i,j} \quad \begin{array}{l}\forall\, i = 2 \ldots I,\\ j = 1 \ldots J_i,\\ c = 1 \ldots n\end{array} \tag{18}$$

$$\theta_c = \sum_{\hat{j}=1}^{J_I}(W_{\hat{j},I+1,1}\sigma_{I,\hat{j},c}) + B_{I+1,1} \quad \forall\, c = 1 \ldots n \tag{19}$$

The objective function (14) incorporates the embeddings outputs $\theta_c$ instead of $X_c^2$, constraints (15) and (16) stay the same as the original model, while constraints (17) through (19) control the neural network embedding for every item class, and are based on the expressions (8) through (10).

The formulation that considers (14) through (19), onward referred as ReLU+, was preliminarily tested against the typical method for embedding a neural network into an optimization model (referred as Classic ReLU), and an off-the-shelf piecewise approximation functionality [9] that comes included in the Gurobi solver 10.0 (referred as gurobiGC); the models for classic ReLU and gurobiGC can be found in the annex 1 and annex 2 respectively.

### Preliminary experiments

Regarding the instances for the preliminary experiments, 4 mayor group has been determined varying the value of $n$ in other to test how the different models sizes affect the computation time, thus the instances were created considering $n \in \{10,100,1000,10000\}$; for each value of $n$, 30 instances were created randomizing the rest of the problem parameters and changing the seed for random number generation. The randomized problem parameters were determined as follow:

- $v_c$: random integer between 50 and 150 distributed uniformly $\forall\, c = 1\ldots n$
- $s_c$: random integer between 10 and 20 distributed uniformly $\forall\, c = 1\ldots n$
- $p_c$: random integer between 5 and 15 distributed uniformly $\forall\, c = 1\ldots n$
- $m_c$: random integer between 2 and 10 distributed uniformly $\forall\, c = 1\ldots n$

It is worth noticing that, because the parameter $m_c$ will be a random integer between 2 and 10, then the values of $X_c$ can only be integers between 0 and 10; this led for the need of training the embedded neural network with those values in mind as inputs, with their corresponding square values as outputs. The weights and biases of the trained neural network can be viewed in anex 3 and 4 respectively.

The experiments were conducted on a computer with a vPRO i5 CPU, that contains 10 cores, each running at 1.6 GHz, with 16 GB of RAM, and a 64bit Windows 11 operating system. The models were programed using Python 3.0 using the library gurobiPy and using an academic free license of GUROBI 10.0 solver. All instances were let to run until optimality was reached in the 3 different models, and the results are shown in table 1

*Table 1: Experiment Results*

|  | | $n =$ | | | |
|---|---|---|---|---|---|
|  | **Type** | **10** | **100** | **1000** | **10000** |
| **Average Time** | gurobiGC | 0.33 | 1.76 | 5.02 | 148.81 |
|  | Classic ReLU | 0.24 | 1.07 | 5.16 | 84.79 |
|  | ReLU+ | 0.13 | 0.45 | 1.59 | 35.32 |
| **Time Standard deviation** | gurobiGC | 0.21 | 1.11 | 1.51 | 79.64 |
|  | Classic ReLU | 0.15 | 0.52 | 1.75 | 21.42 |
|  | ReLU+ | 0.07 | 0.25 | 0.80 | 9.85 |

As it can be seen in table 1, our proposed formulation (ReLU+), is consistently faster than the other models in finding the optimal solutions.

## Conclusions

As this is still a work in process, the preliminary evidence provided and the claims of faster speed that archives the proposed neural network embedding formulation that are made in this technical report should be taken with caution, as further evidence should and will be provided in order to properly sustain said claims. Additionally, the proposed neural network embedding formulation has various limitations that should be taken into account, making it a very niche alternative, and will never be able to fully replace the classical neural network embedding formulation, or the usage of other tools to linearize non linear functions in optimization models. Lastly, current work in progress had shown that the limitation imposed over the proposed neural network embedding weights, that is that $W_{j,i,j} \geq 0$, will make harder the training process of the neural network, especially if the function $G(x)$ is not an increasing function.

**Anex 1**

Below is presented the Classic ReLU formulation; for this formulation the next variables $y_{i,j,c}$ must be added, which is described as follow:

$y_{i,j,c}$: binary variable that takes the value of 1 when $\sigma_{i,j,c} > 0$, and takes the value of 0 when $\sigma_{i,j} \leq 0$; for $i = 1 \dots I, j = 1 \dots J_i, c = 1 \dots n$, and $y_{i,j} \in \{0,1\}$.

$$max \sum_{c=1}^{n}(v_c X_c - p_c(\theta_c - X_c))$$

Subject to:

$$\sum_{c=1}^{n} s_c X_c \leq S$$

$$X_c \leq m_c \qquad \forall\, c = 1 \dots n$$

$$\sigma_{1,j,c} \geq W_{1,1,j} X_c + B_{1,j} \qquad \forall\, j = 1 \dots J_i, \quad c = 1 \dots n$$

$$\sigma_{1,j,c} \leq M y_{1,j,c} \qquad \forall\, j = 1 \dots J_i, \quad c = 1 \dots n$$

$$\sigma_{1,j,c} \leq W_{1,1,j} X_c + B_{1,j} + M(1 - y_{i,j,c}) \qquad \forall\, j = 1 \dots J_i, \quad c = 1 \dots n$$

$$\sigma_{i,j,c} \geq \sum_{\hat{j}=1}^{J_{i-1}} \left(W_{\hat{j},i,j} \sigma_{i-1,\hat{j},c}\right) + B_{i,j} \qquad \forall\, i = 2 \dots I, \quad j = 1 \dots J_i, \quad c = 1 \dots n$$

$$\sigma_{i,j,c} \leq M y_{i,j,c} \qquad \forall\, i = 2 \dots I, \quad j = 1 \dots J_i, \quad c = 1 \dots n$$

$$\sigma_{i,j,c} \leq \sum_{\hat{j}=1}^{J_{i-1}} \left(W_{\hat{j},i,j} \sigma_{i-1,\hat{j},c}\right) + B_{i,j} + M(1 - y_{i,j,c}) \qquad \forall\, i = 2 \dots I, \quad j = 1 \dots J_i, \quad c = 1 \dots n$$

$$\theta_c = \sum_{\hat{j}=1}^{J_I}(W_{\hat{j},I+1,1}\sigma_{I,\hat{j},c}) + B_{I+1,1} \qquad \forall\, c = 1\ldots n$$

## Anex 2

Below is presented the GurobiGC formulation; this formulation is an off the shelf functionality present in gurobi 10.0, in which a picewise linear approximation of a power function is added to the model, further information on this can be viewed on [9] .For the gurobiGC model, a new variable $\mathcal{Y}_c$ must be added , which is presented as follow:

$\mathcal{Y}_c$: value yielded by the piecewise approximation functionality for the observation $c$, where $\mathcal{Y}_c = X_c^2 \ \forall\, c = 1\ldots n$

$$\max \sum_{c=1}^{n}(v_c X_c - p_c(\mathcal{Y}_c - X_c))$$

Subject to:

$$\sum_{c=1}^{n} s_c X_c \leq S$$

$$X_c \leq m_c \qquad \forall\, c = 1\ldots n$$

Its worth noticing that the picewise approximation functionality was defined in python, according to the documentation provided in [9] , as:

modelo.addGenConstrPow(X[c], Y[c], a, "gf"+str(c), "FuncPieces=10")

## Anex 3

In the table below are presented the weights of the neural network used in the embedding of the ReLU+ and Classic ReLU models. It is worth noticing that the network had the following configuration: 1 input value in the input layer(or layer 0), 3 nodes in hidden layer 1, 10 nodes in hidden layer 2, and 1 node in the output layer(or layer 3). Additionally, the information of the table is arranged according to the nomenclature presented for the parameter $W_{\hat{j},i,j}$ , and that all values $W_{\hat{j},i,j} \geq 0$ so ReLU+ works properly.

| $\hat{j}$ | $i$ | $j$ | $W_{\hat{j},i,j}$ |
|---|---|---|---|
| 0 | 1 | 0 | 0.760635 |
| 0 | 1 | 1 | 0.604756 |
| 0 | 1 | 2 | 0.256068 |
| 0 | 2 | 0 | 1.089702 |
| 0 | 2 | 1 | 0.986268 |
| 0 | 2 | 2 | 0.576888 |
| 0 | 2 | 3 | 0.60086 |
| 0 | 2 | 4 | 0.564758 |
| 0 | 2 | 5 | 0.639747 |
| 0 | 2 | 6 | 1.079264 |
| 0 | 2 | 7 | 0.56959 |
| 0 | 2 | 8 | 0.753178 |
| 0 | 2 | 9 | 0.971441 |
| 1 | 2 | 0 | 0.621808 |
| 1 | 2 | 1 | 0.838339 |
| 1 | 2 | 2 | 0.909979 |
| 1 | 2 | 3 | 1.023544 |
| 1 | 2 | 4 | 0.283128 |
| 1 | 2 | 5 | 0.373626 |
| 1 | 2 | 6 | 0.403313 |
| 1 | 2 | 7 | 0.856601 |
| 1 | 2 | 8 | 0.427386 |
| 1 | 2 | 9 | 1.1901 |
| 2 | 2 | 0 | 0.396743 |
| 2 | 2 | 1 | 0.370025 |
| 2 | 2 | 2 | 0.141922 |
| 2 | 2 | 3 | 0.189096 |
| 2 | 2 | 4 | 0.73266 |
| 2 | 2 | 5 | 0.584843 |
| 2 | 2 | 6 | 0.349937 |
| 2 | 2 | 7 | 0.245703 |
| 2 | 2 | 8 | 0.585346 |
| 2 | 2 | 9 | 0.399438 |
| 0 | 3 | 0 | 1.525779 |
| 1 | 3 | 0 | 0.962419 |
| 2 | 3 | 0 | 1.312662 |
| 3 | 3 | 0 | 1.053186 |
| 4 | 3 | 0 | 3.609872 |
| 5 | 3 | 0 | 2.894976 |
| 6 | 3 | 0 | 1.789157 |
| 7 | 3 | 0 | 1.39552 |
| 8 | 3 | 0 | 2.289972 |
| 9 | 3 | 0 | 1.297947 |

**Anex 4**

In the table below are presented the biases of the neural network used in the embedding of the ReLU+ and Classic ReLU models. It is worth noticing that the network had the following configuration: 1 input value in the input layer(or layer 0), 3 nodes in hidden layer 1, 10 nodes in hidden layer 2, and 1 node in the output layer(or layer 3). Additionally, the information of the table is arranged according to the nomenclature presented for the parameter $B_{i,j}$.

| $i$ | $j$ | $B_{i,j}$ |
|---|---|---|
| 1 | 0 | -1.45461 |
| 1 | 1 | -0.39243 |
| 1 | 2 | -2.28163 |
| 2 | 0 | -2.86355 |
| 2 | 1 | 0.000483 |
| 2 | 2 | 0.000182 |
| 2 | 3 | 0.00036 |
| 2 | 4 | -3.78215 |
| 2 | 5 | -3.75593 |
| 2 | 6 | -3.42007 |
| 2 | 7 | 0.000169 |
| 2 | 8 | -3.54801 |
| 2 | 9 | -2.45544 |
| 3 | 0 | 0.141981 |